\input amstex
\documentstyle{amsppt}
\NoRunningHeads
\input pstricks



\magnification=1200

\pagewidth{12.4cm}
\pageheight{18.5cm}
\hcorrection{0.6cm} 
\vcorrection{0.5cm}


\def\today{\ifcase\month\or 
 January\or February\or March\or April\or May\or June\or 
 July\or August\or September\or October\or November\or 
December\fi 
 \space\number\day, \number\year} 


\define\newlabel#1#2{\expandafter\ifx\csname MyLabel#1\endcsname\relax\else
 \immediate\write16{Label `#1' is multiply defined.}\fi
 \global\expandafter\def\csname MyLabel#1\endcsname{#2}}

\newif\ifprintlabels
\printlabelsfalse

\define\Ref#1{\expandafter\ifx\csname MyLabel#1\endcsname\relax
 \immediate\write16{Reference `#1' is undefined.}{\bold ??}\else
 \csname MyLabel#1\endcsname\ifprintlabels\smash{\rlap{${}^*$}}\fi
 \fi}

\define\eqref#1{(\Ref{#1})}

\define\biblabel#1{\ifprintlabels\smash{\llap{\fiverm #1\ }}\fi\Ref{#1}}

\define\label#1{{\ifprintlabels\smash{\raise0.8\baselineskip
 \rlap{\fiverm #1}}\fi\Ref{#1}}}

\topmatter
\title  Tits geometry, arithmetic groups and the proof of a conjecture of Siegel
\endtitle
\author Enrico Leuzinger \endauthor
\address  Mathematisches Institut II,  Universit\"at Karlsruhe (TH), 
Englerstrasse 2, D--76131 Karlsruhe, Germany \endaddress
\email Enrico.Leuzinger\@math.uni-karlsruhe.de\endemail
\input amstex
\documentstyle{amsppt}
\NoRunningHeads
\input pstricks



\magnification=1200

\pagewidth{12.4cm}
\pageheight{18.3cm}
\hcorrection{0.6cm} 
\vcorrection{0.5cm}


\def\today{\ifcase\month\or 
 January\or February\or March\or April\or May\or June\or 
 July\or August\or September\or October\or November\or 
December\fi 
 \space\number\day, \number\year}


\abstract  Let $X = G/K$ be a Riemannian symmetric space of noncompact type and of rank $\geq 2$. 
 An irreducible, non-uniform  lattice $\Gamma\subset G$ in the isometry group 
of $X$ is arithmetic and gives rise to a locally symmetric space $V=\Gamma\backslash X$.
Let  $\pi:X\rightarrow V$ be the canonical projection.
Reduction theory for arithmetic groups provides a dissection $V=\coprod_{i=1}^k \pi(X_i)$ 
with 
$\pi(X_0)$ compact and such that $\pi$ restricted to $X_i$ for $i\geq 1$ is injective.
In this paper we complete reduction theory by focusing on metric properties of the sets 
$X_i$. We detect subsets $C_i$ of $X_i$  (${\Bbb Q}$--Weyl chambers) such that $\pi_{\mid C_i}$ 
is an isometry and such that $C_i$ is a net in $X_i$. This result is then used to prove a conjecture 
of C.L. Siegel. We also show that the locally symmetric space $V$ is quasi-isometric to 
the Euclidean cone over a finite simplicial complex  and study the Tits geometry of $V$.
\endabstract
\endtopmatter

\document

\head  1. Introduction and survey of results\endhead

In this paper we study geometric properties of certain locally symmetric spaces 
$V = \Gamma\backslash X$. Here $X$ is  a Riemannian symmetric space 
of
noncompact type (without Euclidean factor)  and  $\Gamma$
 is a  non-uniform (torsion-free)  lattice  in the group of isometries of $X$, i.e., the  
quotient space $V$ is not compact but has finite volume. We  make the additional assumption 
that the rank of $X$ is at least two and we also
assume  that the non-uniform lattice $\Gamma$ is irreducible. In that case a
fundamental theorem of G.A. Margulis  asserts that  $\Gamma$ is {\it arithmetic}. Roughly 
this means that there is a semisimple
 algebraic group 
${\bold G}\subset \bold{GL}(n, \Bbb C)$ defined over $\Bbb Q$ such that $\Gamma$ is 
commensurable with
${\bold G}(\Bbb Z) = {\bold G}\cap  \bold{GL}(n, \Bbb Z)$ (see [Z] Ch. 6 and Section 2 below).
The symmetric space $X$  can be written as
the homogeneous space ${\bold G}({\Bbb R})^0/K$ for a (chosen) maximal compact 
subgroup $K$. The ${\bold G}({\Bbb R})^0$-invariant Riemannian metric on $X$ induces one on $V$ so
 that the canonical projection
$\pi:X\longrightarrow \Gamma\backslash X$ is a Riemannian covering.

\subhead 1.1.  Metric properties of locally symmetric spaces \endsubhead

Fundamental domains provide a basic and classical tool to investigate quotients of
discrete groups acting isometrically on spaces of constant curvature.  For lattices in 
isometry groups of  higher rank  symmetric spaces {\it exact} fundamental domains have a 
rather  complicated geometry. 
In order to study  locally symmetric spaces $V = \Gamma\backslash X$ one therefore
uses only  coarse  fundamental domains, whose shapes reflect  certain (algebraic) 
decompositions of the isometry group of $X$.

{\it Reduction theory} provides such
coarse fundamental domains $\Omega\subset X$ for arithmetic groups. Its basic building blocks 
are so-called Siegel sets, which in turn are unions of certain Euclidean cones (
${\Bbb Q}$--Weyl chambers). A more refined version, 
{\it precise reduction theory}, yields a combinatorial dissection of $V$ into disjoint subset:
$V=\coprod_{i=0}^l V_i$ where $V_0$ is compact and the other $V_i$ are noncompact
and can be explicitely described as one-to-one images of (generalized) Siegel sets, $V_i=\pi({\Cal S}_i)$,
where  $\pi:X\longrightarrow  V=\Gamma\backslash X$ denotes the canonical 
projection map (see e.g. [S], 9.6). A variant of such a dissection based on an
 exhaustion of $V=\Gamma\backslash X$ by compact 
submanifolds with corners has been created in [L2] and in [L3] and  will be our main tool
in investigating $V$. We shall discuss reduction theory and exhaustions in more detail
in Sections 2 and 3.

Our  main result (Theorem 4.1) completes precise reduction theory by revealing
 {\it metric properties} of subsets of the  above  $V_j$: 

\proclaim{Theorem A}  There are ${\Bbb Q}$--Weyl chambers ${\Cal C_i}\subset {\Cal S}_i, 
1\leq i\leq l$, in $X$ such that,
for each $i$,  the restriction of the $\pi:X\longrightarrow V$ to ${\Cal C_i}$ 
is an  isometry (with respect to the distance functions 
induced by the Riemannian metrics on $X$ resp. $V$). Moreover
 $\bigcup_{i=1}^l\pi({\Cal C_i})$ is a net in $V$.
\endproclaim

As a consequence of Theorem A we
confirm a conjecture of C.L.Siegel.

\proclaim{Theorem B} The projection map $\pi:X\longrightarrow \Gamma\backslash X$ restricted
 to a
Siegel set is a $(1,D)-$quasi-isometry.
\endproclaim

 Siegel's conjecture has been proved independently by 
L.\,Ji  and in the special case of $SL(n,{\Bbb Z})\subset SL(n,{\Bbb R})$ also by J.\,Ding 
(see [J] and [D]).
It has been applied recently by D.\,Kleinbock and G.A.\,Margulis to investigate logarithm laws for 
 certain flows on 
locally symmetric spaces (see [KM]).

\subhead 1.2. Quasi-isometry invariants for  locally symmetric spaces\endsubhead

Using the assertion of  Theorem A  we  construct a quasi-isometric model  of the locally symmetric 
space $V=\Gamma \backslash X$. It is obtained as the  {\it 
asymptotic cone} which is defined as the Gromov-Hausdorff limit of ``rescaled'' pointed 
metric spaces:
$$ \text{Cone}(V) :=  {\Cal H}-{\lim}_{n\rightarrow\infty} (V, v_0, \frac{1}{n} d_V),$$
where $v_0$ is an arbitrarily chosen  point of $V$ and $d_V$ is the distance function on 
$V$ induced by the Riemannian metric (see e.g. [GLP] Ch. 3 or [Bu] Ch. I).  We remark that 
in contrast to the case considered here the   definition of an asymptotic
 cone in general involves the use of ultrafilters, and the limit space may even  depend on the 
chosen 
ultrafilter.
 Various aspects of  asymptotic cones of general spaces are discussed in  M.\,Gromov's essay [G].
 Recently 
B.\,Kleiner and B.\,Leeb also used  
 asymptotic cones to prove the rigidity of quasi-isometries of symmetric spaces  (see [KL]).

In some cases the asymptotic cone of $\Gamma\backslash X$ is  easy to describe. For 
example 
if the rank of $X$ is 1 (or, more generally, if $V$ is a Riemannian manifold with finite 
volume and strictly negative sectional curvature),  then   $\roman{Cone}(V)$ is 
a ``cone''  over $k$ points, i.e., $k$ rays with a common origin, where $k$ is 
the number of ends of $V$.  For  a Riemannian product $V= V_1\times V_2$, where 
$V_1$, $V_2$ are as in the previous example and each has only one end,   
$\roman{Cone}(V)$  can be identified with the first quadrant in $\Bbb R^2$.
Much more intricate are quotients of $SL(n,\Bbb R)/SO(n)$ by congruence subgroups of 
$SL_n(\Bbb Z)$. These are  examples of higher rank locally symmetric spaces which are not 
only products of rank $1$ factors; their asymptotic geometry was first studied by T.\,Hattori
(see [H1], [H2]). In what follows we are concerned with the general higher rank case.

By the Cartan-Hadamard theorem a simply connected, complete Riemannian manifold $X$ of
 nonpositive sectional curvature 
 is diffeomorphic to ${\Bbb R}^n$,
$n = \dim X$. One can  compactify $X$  by adjoining the  {\it boundary at 
infinity},
$\partial_{\infty}X\cong S^{n-1}$, which is defined as the set of equivalence classes of asymptotic 
geodesic rays. The enlarged space $X\sqcup
 \partial_{\infty}X$ is homeomorphic to a closed ball in ${\Bbb R}^n$. 
On  $\partial_{\infty}X$ there is a natural distance function, the {\it Tits metric}, which 
is
also defined in terms of geodesic rays (see Section 3). 
In the special case where $X$ is a symmetric space of noncompact type the 
boundary at infinity $\partial_{\infty}X$ carries the additional structure of a 
spherical Tits bulding given by the disjoint union of Weyl chambers and their walls ``at 
infinity''. In that  case the building structure  on $\partial_{\infty}X$ can be 
reconstructed using  the Tits metric (see [BGS]).

Since the algebraic group $\bold G$ is defined over $\Bbb Q$ one also 
has the {\it Tits building $\Cal T$ of 
$\bold G$ over} $\Bbb Q$, which  is
a simplical complex whose combinatorial structure reflects the incidence relations
among parabolic ${\Bbb Q}$--subgroups of $\bold G$. The lattice $\Gamma$ acts naturally on   
a geometric realization
 $|\Cal T|$ of $\Cal T$ as a spherical simplicial complex. The resulting quotient 
$\Gamma\backslash |\Cal T|$ turns out to be a finite  simplical complex. 
If $\text{rank}_{ \Bbb Q 
}\Gamma \geq 2$ the spherical complex $|\Cal T|$ is  equipped with a spherical  
distance function  which we normalize such that the diameter of $|\Cal T|$ is 
$\pi$  (see Section 6). This  metric in turn  induces a distance function 
$d_{\Cal T}$ on the quotient space  $\Gamma\backslash |\Cal T|$.  The  {\it 
Euclidean cone}  $C(\Gamma\backslash |\Cal T|)$  {\it over }
 $\Gamma\backslash |\Cal T|$ is the product $[0,\infty)\times  \Gamma\backslash 
|\Cal T|$ with $\{0\}\times \Gamma\backslash |\Cal T|$ collapsed to a point and  
endowed with the  cone metric
$$d_C^2((a,x),(b,y)) := a^2 + b^2 - 2ab\cos d_{\Cal T}(x,y).$$

We remark that  $\Gamma\backslash |\Cal T|$ is connected if 
 $\Bbb Q$--rank of $\Gamma$ is $\geq 2$.

If  $\text{rank}_{ \Bbb Q }\Gamma =1$,  $|\Cal T|$ consists of countably many 
points and $C(\Gamma\backslash |\Cal T|)$ is a cone over a finite number of 
points as in the first example described above. 
In that  case we set $d_C((a,x),(b,y)) := a+b$ if $x\neq y$ and 
$d_C((a,x),(b,y)):=|a-b|$ if $x=y$. 
 
We can now describe  the asymptotic  cone of an arithmetic quotient:

\proclaim{Theorem C} Let $X$ be a symmetric space of noncompact type (without Euclidean 
factor) and rank $\geq 2$. Let $\Gamma$ be an irreducible, non-uniform lattice in the group 
of isometries of $X$. Then the
 asymptotic cone of the locally symmetric space $V=\Gamma\backslash X$ is quasi-isometric to $V$ and 
 is isometric to 
the Euclidean cone over the finite simplicial complex $\Gamma\backslash |\Cal T|$.
\endproclaim

Various invariants of the algebraic group $\bold G$ are reflected in the geometry of
 $X = {\bold G}(\Bbb R)^0/K$ and of  $\Gamma\backslash X$. For instance it
 is well-known that the $\Bbb R$--rank of $\bold G$  equals the rank of the symmetric space 
$X$,
i.e., the maximal dimension of a totally geodesic, flat submanifold  in $X$.

 It is an elementary but  crucial fact that
bi-lipschitz invariants of  the asymptotic cone of $V=\Gamma\backslash X$ are quasi-isometry 
invariants of $V$.
Theorem C  thus yields the following geometric interpretation of the (purely algebraic) 
notion of ``$\Bbb Q$--rank'' as a quasi-isometry invariant of $V$. Note that it 
follows from  our assumptions  that the $\Bbb Q$--structure of $G$ is unique (see 
[T1]) and we can thus set $\text{rank}_{ \Bbb Q }\Gamma :=
\text{rank}_{ \Bbb Q} {\bold G}$.  

\proclaim{Corollary D}
Let $\Gamma$ be an arithmetic lattice as in Theorem C. Then 
$$\roman{rank}_{\Bbb Q} \Gamma = \dim \, \roman{Cone} (\Gamma\backslash X).$$
\endproclaim

The dimension of the Euclidean cone  $C(\backslash |\Cal T|)$ is defined  as the 
(topological) dimension of the simplicial complex $\Gamma\backslash |\Cal T|$ {\it plus 
one}. According to [KL] the asymptotic cone(s) of a 
globally symmetric space $X$  are Euclidean buildings and 
 one has $\dim  \roman{Cone} (X) = \roman{rank}_{\Bbb R} {\bold G}$.

The asymptotic cone of a compact metric space is just a  point. Corollary D   thus  reflects 
geometrically  a classical 
 result of A.\,Borel and Harish-Chandra asserting that 
an arithmetic lattice  $\Gamma$ is co-compact, i.e., $\Gamma\backslash X$ is compact, if and 
only if  the $\Bbb Q$--{\rm rank} of $\Gamma$ {\rm (}or $\bold G${\rm )} is $0$ (see [BH]).

Typical for Gromov-Hausdorff limits are collapsing phenomena (as for instance the  just 
mentioned fact that the  asymptotic cone of a compact space is a  point). This makes the 
limit objects in some sense easier to understand. On the other hand  
Gromov-Hausdorff limits of manifolds usually degenerate into (singular)
metric spaces. In particular concepts which are related to the differentiable structure, 
like
sectional curvature, may have no meaning in the limit. It is therefore 
necessary to work in a more general setting using notions that are  defined for 
(geodesic) metric spaces. ``Nonpositive sectional curvature'', for instance, 
can be generalized by the local property  ``curvature bounded from above by $0$'' (resp. by its 
global counterpart,  ``CAT($0$) space''). Roughly speaking one requires that sufficiently small 
(resp. all) 
geodesic triangles
are  ``thinner'' than their comparison triangles of the same side length in the Euclidean 
plane. For 
example a locally symmetric space of non compact type has curvature bounded 
above by $0$  and every  globally symmetric space is  
CAT($0$).  Recent references for these purely metric
concepts   due to A.D.\,Alexandrov are [Ba] and the comprehensive [BrHa].

 Notice that $V=\Gamma\backslash X$ has curvature bounded above by $0$, but that $V$ (in 
contrast to $X$) is not a CAT($0$) space since $V$ is in general not contractible. However 
we have the following corollary to Theorem C.

\proclaim{Corollary E} The asymptotic cone of a locally symmetric space $\Gamma\backslash X$ 
 (as in Theorem C) is a $\roman{CAT(0)}$ space. 
\endproclaim

The  boundary at infinity $\partial_{\infty}V$ of $V = \Gamma\backslash X$ is 
defined as the set of 
equivalence classes of asymptotic geodesic rays in $V$.
From Theorem C we deduce  the following combinatorial description of the boundary 
at infinity of a locally symmetric space, which  was first proved by L.\,Ji and 
R.\,MacPherson (see [JM] and also [H2], [L1]):

\proclaim{Theorem F}
Let $V = \Gamma\backslash X$ be as in Theorem C. The points of the  boundary at infinity
$\partial_{\infty}V$  correspond bijectively to the points  of 
$\Gamma\backslash |\Cal T|$, i.e.,
 $$\partial_{\infty}V\cong\Gamma\backslash |\Cal T|\cong\partial_{\infty}C(\Gamma\backslash 
|\Cal T|).$$
\endproclaim

 The next theorem asserts that the distance  $d_{\Cal T}$ on $\Gamma\backslash |\Cal T|$ can 
be defined intrinsically in terms of geodesic rays in $V$. This answers a 
question of  Ji and MacPherson (see [JM],  Remark 6.4). 

\proclaim{Theorem G} Let $\text{rank}_{ \Bbb Q }\Gamma \geq 2$. For any two points $z_1,z_2$ 
in the  boundary at infinity $\partial_{\infty}V\cong \Gamma\backslash |\Cal T|$ represented 
by unit speed  geodesic rays $c_1, c_2$ in $V$, respectively, one has the formula
$$2\sin\frac{1}{2} d_{\Cal T}(z_1,z_2) = 
\lim_{t\rightarrow\infty}\frac{1}{t}d_V(c_1(t),c_2(t)).$$
\endproclaim

Observe that  this is exactly the same formula as the one for  simply connected, 
complete Riemannian   manifolds of nonpositive curvature or more generally for  CAT($0$) 
spaces (see [BGS] Ch. 4, [BrHa] and Section 5 below). For symmetric spaces  the (crucial) 
difference is that  
 in the simply connected case $d_{\Cal T}$ refers to the Tits distance on the 
spherical building associated to $\bold G (\Bbb R)$, whereas in the {\it non} simply 
connected
case $d_{\Cal T}$ is induced by the Tits distance of the spherical building associated to 
$\bold G (\Bbb Q)$. While the former (beeing independent of any $\Bbb 
Q$--structure on the isometry group of
$X$)  reflects geometric properties of $X$,  the latter depends on 
the $\Bbb Q$--structure  associated to the lattice $\Gamma$ and encodes (large 
scale) geometric properties 
of the quotient $V = \Gamma \backslash X$.

\head 2. Coarse fundamental domains for arithmetic groups\endhead

Let   $G$  denote the
 identity component of the  group of isometries of the symmetric space $X$; it
 is a connected, semisimple Lie group with trivial center.  In this section we
 discuss some classical results of reduction theory about coarse fundamental domains
for arithmetic groups, which are adopted to the group structure of $G$.

We make two natural assumptions
about the non-uniform lattice $\Gamma$ in $G$. Firstly, we assume that $\Gamma$ is   {\it
irreducible}
 (see
[R2] 5.20) and, secondly, that $\Gamma$ is {\it neat} and hence in particular torsion-free
 (see [B1] \S 17). Under these assumptions
  the arithmeticity theorem
of Margulis asserts that there is a connected
semisimple linear algebraic group
 ${\bold G}$ defined over $\Bbb Q$, $\Bbb Q$--embedded in a general linear group
$\bold{GL}(n,\Bbb C)$, and a Lie group isomorphism
$$\rho: G
\longrightarrow {\bold G}(\Bbb R)^0  = (\bold G
\cap \bold {GL} (n,\Bbb R))^0
$$
 such that
$\rho(\Gamma)$ is  {\it arithmetic}. This means that
  $$\rho(\Gamma) \subset {\bold G}(\Bbb Q)  = \bold G
\cap \bold {GL} (n,\Bbb Q)\subset \bold {GL}(n,\Bbb C)$$
 is commensurable
with the group  ${\bold G}(\Bbb Z) = {\bold G}\cap \bold {GL}(n,\Bbb Z)$ (see
[Z] 3.1.6 and 6.1.10). Note that the
 symmetric space $X$ can be recovered as the manifold of maximal compact
 subgroups
of the identity component of the group $\bold  G(\Bbb R)$
 of $\Bbb R$--rational  points of $\bold G$.   For simplicity we
identify $G$ with
${\bold G}(\Bbb R)^0$ and  $\Gamma$ with $\rho(\Gamma)$ in this paper.

\remark{\bf Remark 2.1} For an  irreducible, non-uniform  lattice $\Gamma$ as above the
connected semi\-simple $\Bbb Q$--group $\bold G$ is unique up to isomorphism (see [T1] Section
3.2). We can thus define $\text{rank}_{ \Bbb Q }\Gamma := \text{rank}_{ \Bbb Q} {\bold G}$.
\endremark

\subhead 2.1. Horocyclic coordinates and ${\Bbb  Q}$--Weyl chambers\endsubhead

We recall some definitions and basic facts
 about linear algebraic groups together with geometric interpretations which  are 
needed below.

Let $\bold S$ (resp.
$\bold T$) be
 a maximal $\Bbb Q$--split
(resp.  a maximal $\Bbb R$--split) algebraic torus  of
  $\bold G$, i.e. a subgroup of $\bold G$
which is isomorphic
over $\Bbb Q$ (resp. $\Bbb R$)  to the direct product of $q$
(resp. $r\geq q$) copies of $\Bbb C^*$.
All such tori are conjugate under ${\bold G}(\Bbb Q)$ (resp. ${\bold G}(\Bbb R)$)
and their common dimension  $q$ (resp. $r$) is
called the $\Bbb Q$--{\it rank} (resp. $\Bbb R$--{\it rank}) of
${\bold G}$.
 The identity component
of ${\bold S}(\Bbb R)$ (resp. ${\bold T}(\Bbb R)$)
 will be denoted by $A$
(resp. by $A_0$), the corresponding  Lie algebra by $\frak a$ (resp.
by  $\frak a_0$). The $\Bbb R$--rank of ${\bold G}$ coincides with the
 rank of the symmetric  space
 $X$, i.e., the maximal dimension of  totally geodesic flat subspaces.
The choice of a maximal compact subgroup $K$ of $G$
 is equivalent to the choice of a 
base point $x_0$ of  $X$. 
We can further choose $K$ 
with Lie algebra $\frak k$ so that under the corresponding Cartan decomposition
 $\frak g = \frak k\oplus \frak p$ 
of the Lie algebra $\frak g$ of $G$  we have $\frak a
\subseteq 
\frak a_0\subset 
\frak p\cong T_{x_0}X$. Here
 $\frak a_0$ is  maximal abelian in $\frak p$, i.e., the tangent space
 at $x_0$
 of the (maximal $\Bbb R$--) flat $A_0\cdot x_0$ in $X$.
The pair of Lie algebras $(\frak g, \frak a_0)$
gives rise to  the root system $_{\Bbb R}
\Phi$ of the
symmetric space. Similarly there is a
system of $\Bbb Q$--roots $_{\Bbb Q}\Phi$  associated to
 the pair
$(\frak g, \frak a)$ (see [B3] \S 21).
It is always possible to choose
orderings of $_{\Bbb Q}\Phi$ and $_{\Bbb R}\Phi$ such 
 that the restrictions of  simple  $\Bbb R$--roots 
of $_{\Bbb R}\Phi$ to $\frak a$
are either {\it simple} $\Bbb Q$--roots of $_{\Bbb Q}\Phi$, i.e., the
elements of a basis $\Delta =\  _\Bbb Q\Delta$ of $_{\Bbb Q}\Phi$,  or zero
 (see [BT] 6.8).
 The basis $_{\Bbb R}\Delta$ defines a closed $\Bbb R$--Weyl chamber 
 $\overline{\frak a_0^+}$  in 
$\frak a_0$ and  $\Delta$ then  determines a
closed {\it $\Bbb Q$--Weyl chamber}  in $\frak a$ 
$$\overline{\frak a^+} 
:= \{H\in \frak a\mid \alpha(H)\geq 0, \
\roman{for\  all}\ \  \alpha\in \Delta\}.
$$
 We set $\overline{A^+} = \exp \overline{\frak a^+}$ (resp.
$\overline{A_0^+} = \exp \overline{\frak a_0^+}$).
 A {\it  $\Bbb Q$--Weyl chamber  in  $X$} is  a  translate of the basic chamber
$\overline{A^+}\cdot x_0\subseteq \overline{A_0^+}\cdot x_0$ under $\bold G(\Bbb Q)$.
 The elements of 
 $\Delta$ are differentials of characters (defined over
 $\Bbb Q$) 
of the
maximal $\Bbb Q$--split torus $\bold S$. It is convenient to
identify the elements of $\Delta$ also  with such characters.  When 
restricted to $A$ we denote their
values  by  $\alpha(a)$ ($a\in A, \alpha 
\in \Delta$).
Notice that $\overline{A^+} = \{a\in A\mid
\alpha(a) \geq 1\ \ \roman{for \ all} \ \ \alpha\in\  \Delta\}$.

A closed subgroup $\bold P$  of $\bold G$ defined over $\Bbb Q$
is a {\it parabolic $\Bbb Q$--subgroup} if $\bold G/\bold P$ is a projective
variety (see  [B3] \S 11). A {\it parabolic $\Bbb Q$--subgroup 
$P$ of} $G = \bold G(\Bbb R)^0$ is by definition the intersection
of $G$ with a parabolic $\Bbb Q$--subgroup of $\bold G$ (see [BS]).
The conjugacy classes under $\bold G$($\Bbb Q$) of parabolic $\Bbb Q$--subgroups
are in one-to-one correspondence with the subsets $\Theta$ of the 
(chosen) set  $\Delta$ of simple $\Bbb Q$--roots; they are represented by 
the {\it standard parabolic} $\Bbb Q$--subgroups $\bold P_{\Theta}$
of $\bold G$ (see [B3]  21.11). The corresponding standard parabolic 
 $\Bbb Q$--subgroups of $G$
are denoted by $P_{\Theta}$. 
Let $\Cal P$ be the collection of all parabolic subgroups of $\bold G$ defined over $\Bbb 
Q$. For 
every $\bold P\in \Cal P$ the corresponding parabolic $\Bbb Q$--subgroup $P = G\cap \bold P$ 
of $G$ has a 
Langlands-decomposition $P= U_{\bold P}M_{\bold P}A_{\bold P}$, where $U_{\bold P}$ is 
unipotent and $M_{\bold P}$ is reductive. Furthermore
$A_{\bold P}$ centralizes $M_{\bold P}$ and normalizes $U_{\bold P}$ (see [BS] and [B1]).
 This yields a (generalized)
Iwasawa decomposition for $G$, i.e., 
 $G = P K = U_{\bold P} M_{\bold P} A_{\bold P}  K$, which  in turn implies 
 that $P$  acts transitively 
on the symmetric space $X$.
The intersection of the maximal compact subgroup $K$ of $G$ with $M_{\bold P}$ is 
maximal compact in $M_{\bold P}$ and the quotient $Z_{\bold P} = M_{\bold P}/(K\cap M_{\bold 
P})$
is in general the Riemannian product of a symmetric space 
of noncompact type by a flat Euclidean space. 
Let $\tau_{\bold P} : M_{\bold P}\longrightarrow Z_{\bold P} $ be the natural projection. 
Then the
``horocyclic coordinate map" 
$$\mu_{\bold P} : Y_{\bold P} = U_{\bold P}\times Z_{\bold P} \times A_{\bold P} 
\longrightarrow X \ \ ;
\ \  (u, \tau_{\bold P}(m), a)\longmapsto u m a\cdot x_0 $$
is an isomorphism of analytic manifolds (see [BS]).

\subhead 2.2. Reduction theory  \endsubhead

A subset $\Omega\subset X$ is called a {\it fundamental set} for an arithmetic group $\Gamma$
if the following two conditions are satified: (i) $X=\Gamma\cdot \Omega$; (ii) for every
 $q\in {\bold G}({\Bbb Q})$ the set $\{\gamma\in \Gamma\mid q\Omega\cap\gamma\Omega\neq\emptyset\}$ is 
finite. It is the goal of {\it reduction theory} to provide explicit fundamental sets.

Let ${\bold P}$ be the  {\it minimal} parabolic $\Bbb Q$--subgroup of ${\bold G}$
which corresponds to the  $\Bbb Q$--Weyl chamber $\overline{A^+}\cdot x_0$.
 A {\it
generalized Siegel set} $\Cal S =  \Cal S_{\omega,\tau}$
  in $X$ (relative to  $\overline{A^+}\cdot x_0$)
is a subset of $X$ of the form
$\Cal S_{\omega, \tau} = \omega A_{\tau}\cdot x_0$ where $\omega$ is a
relatively compact in $U_{\bold P} M_{\bold P}$ 
and, for $\tau>0$,  $A_{\tau} = \{a\in A = A_{\bold P}\mid\alpha(a)\geq \tau,\ \alpha\in 
\Delta\}$.
If we define $a_0\in A$ by $\alpha(a_0) = \tau$ for all
$\alpha\in \Delta$, then $A_{\tau} = A_1a_0 = \overline{A^+}a_0$
and ${\Cal C} = A_{\tau}\cdot x_0\subset \Cal S$ is a (translate of a) 
$\Bbb Q$--Weyl chamber in $X$.
   Fundamental sets for arithmetic
groups consist of finitely many Siegel sets (see  [B1] \S 13 and \S 15 for the proofs):

\proclaim{Proposition  2.2} ({\rm Borel, Harish-Chandra}) Let $\bold G$ be a
semisimple  algebraic group defined over $\Bbb Q$ with associated
Riemannian symmetric
 space $X = G/K$ (of maximal compact subgroups of $G$). Let $\bold P$ be a
minimal parabolic $\Bbb Q$--subgroup of $\bold G$ and let $\Gamma$ be an 
arithmetic subgroup of $\bold G(\Bbb Q)$. Then
  there exists a generalized Siegel set
 $\Cal S = \Cal S_{\omega,\tau}$
 such that,  for a fixed set 
$\{q_i\mid 1\leq i\leq m\}$ of representatives of the  finite  set of double cosets 
$\Gamma\backslash 
 \bold G(\Bbb Q)/\bold P(\Bbb Q)$, the union
$\Omega = \bigcup_{i = 1}^m q_i\cdot{\Cal S}$ is a fundamental set
 of finite volume for 
$\Gamma$ in $X$.
\endproclaim

\head 3. Polyhedral exhaustions of locally symmetric spaces \endhead

A crucial ingredient in the proof of 
the main (technical) Theorem 4.1 below is the fact that a  locally symmetric space 
$V = \Gamma\backslash X$ can  be exhausted (in many ways) by  polyhedra, i.e., by  compact 
submanifolds $V(s)$ of $V$ with corners: $V = \bigcup_{s\geq 0}V(s)$. 
Such exhaustions were constructed in [L2] and further investigated in [L3].
Similar constructions appear also in papers by J.\, Arthur, M.S.\, Osborne-G.\, Warner and 
L.\, Saper (see  [A], [OW] and [S]).

In this section we recall some facts about the geometric structure of the polyhedra 
$V(s)$ for $s\geq 0$. We invite the reader to consult the articles 
[L2] and [L3] for more details and proofs.

\subhead 3.1. Exhaustions and the rational Tits building  \endsubhead

 Let $d_X$ be the  distance function  induced by the 
Riemannian metric of the symmetric space $X$. A unit-speed {\it geodesic ray} is a curve  
$c : [0,\infty)\rightarrow X$ 
which realizes the distance 
between any two of its points, i.e., $c$ is an isometric embedding of $[0,\infty)\subset 
\Bbb R$ into $X$. We call two unit-speed  geodesic rays $c_1, c_2 :[0,\infty)\rightarrow X$ 
in $X$ {\it asymptotic} if there is a constant $C\geq 0$ such that 
$\lim_{t\rightarrow\infty}d_X(c_1(t),c_2(t))\leq C$; the corresponding equivalence class is
denoted by $c_1(\infty) (= c_2(\infty))$. The {\it  boundary of $X$ at
infinity}, $\partial_{\infty}X$, is  defined as the set of equivalence classes of asymptotic 
geodesic rays. 
A {\it Busemann function} $h$ associated to a  ray $c : [0,\infty)\rightarrow X$ 
representing
 the equivalence class $z = c(\infty) \in \partial_{\infty}X$ is defined by
$$h : X\rightarrow \Bbb R\ ;\ \  x\longmapsto \lim_{t\rightarrow\infty}[ d_X(c(t),x) - t].$$
A {\it horosphere} (resp. a {\it horoball}) centered at $z\in \partial _{\infty}X$  is a 
level set (resp. sublevel set) of $h$.

We next recall the definition of the (rational) Tits building $\Cal T$ of $\bold G$.
(see e.g. [T2] 5.2).
Let   ${\Cal P}_m$ be the set 
of all  {\it maximal} parabolic $\Bbb Q$--subgroups of $\bold G$ (and $\neq \bold G$) (see 
Section 2).
The {\it Tits building $\Cal T$ associated to $\bold G$ over $\Bbb Q$} is the simplicial
complex whose set of vertices is ${\Cal P}_m$, and whose simplices are the non-empty subsets
$I$ of ${\Cal P}_m$ such that $\bold P_I = \bigcap_{\bold P\in I}\bold P$ is a parabolic
$\Bbb Q$--subgroup of $\bold G$. A geometric realization  
$|\Cal T|$ of $\Cal T$ in  $\partial_{\infty} X$ is obtained as follows. We define a 
(closed) {\it $\Bbb Q$--Weyl
chamber at infinity} as
$$\split (\overline{A^+} x_0)(\infty) & := \{ c_H(\infty)\mid c_H(t) = \exp tH\cdot x_0, \ 
H\in
 \overline{\frak a^+}, \|H\| =1 \}\cong \\
&\cong\{H\in \overline{\frak a^+}\mid \|H\| =1 \}.
\endsplit
$$
Any other $\Bbb Q$--Weyl chamber at infinity is of the form $g\cdot (\overline{A^+} 
x_0)(\infty)$ with $g\in \bold G (\Bbb Q)$. The geometric realization $|\Cal T|$
  is then defined as the union of all $\Bbb Q$--Weyl chambers and their walls
 in $\partial_{\infty}X$ partially ordered with the order relation ``wall of'' (see also 
[JM]).

The following proposition is proved in [L3], Theorem 3.6.

\proclaim{Proposition 3.1} Let $V = \bigcup _{s\geq 0}V(s)$ be a polyhedral  
exhaustion as introduced above. For any $s\geq 0$, the polyhedron $V(s)$ can be written as 
$\Gamma\backslash X(s)$,  where $X(s)$ is the complement  of a $\Gamma$--invariant union of
 countably many open horoballs: 
$$
X(s) = X - \bigcup_{k}{\Cal B}_{k}(s).
$$
The ``centers'' in $\partial_{\infty}X$  of the deleted horoballs ${\Cal 
B}_{k}(s)$ are in one-to-one  correspondence with the maximal parabolic 
subgroups of ${\bold G}$ defined over $\Bbb Q$, or, equivalently, with  the 
vertices of a geometric realization of the Tits building $\Cal T$ of $\bold G$.
\endproclaim

If the $\Bbb Q$--rank of $\Gamma$ is $1$, then the deleted horoballs ${\Cal B}_{k}(s)$ are 
disjoint in $X$. If
 the $\Bbb Q$--rank of $\Gamma$ is $\geq 2$, then, in contrast, the  horoballs ${\Cal 
B}_{k}(s)$ intersect and
give rise to the corners of $X(s)$ (and $V(s)$). 
Their local structure  can  be described  for example by taking the (inner) unit normals of the 
involved  horospheres $\partial {\Cal B}_{k}(s)$ (see also [L2] Lemma 4.1 (iv)).
In fact, any of these horospheres  is of the form
$\gamma\{\tau_{\alpha}^{-1}\tilde{h}_{i\alpha} = -s\}$
with $\gamma\in \Gamma$, $i\in\{1,\ldots,m\}$, $\alpha\in\Delta$ and where the 
$\tilde{h}_{i\alpha}$ are  Busemann functions on $X$ associated to distinguished  geodesic 
rays in the fundamental set $\Omega\subset X$  (see [L2], Theorem 3.6).

The inner unit normal field of the horosphere
 $ \{\tau_{\alpha}^{-1}\tilde{h}_{i\alpha} = -s\}$ in $X$ is given by
$ -\roman{grad}\ \tilde{h}_{i\alpha}$ (see  [HI] Proposition 3.1).
Let  $p$ be a point of $V(s)$.
The {\it outer angle} $O(p)$
 at  $p$  is defined as the set  of all unit tangent vectors $v\in T_p V(s)$ such that
$\langle v, w\rangle_p\leq 0$ for all $w$ in the tangent cone of $V(s)$ at $p$.
Let $\Theta$ be a subset of $\Delta$.
If the horospheres $\{\tau_{\alpha}^{-1}\tilde{h}_{i\alpha} = -s\}$, for $\alpha\in \Theta$, 
intersect, then the outer angle $O(p)$ at an intersection point $p$
 can be identified with all  positive linear combinations
(of norm $1$) of the
unit vectors 
$-\roman{grad}\ \tilde{h}_{i\alpha}(p)$,  $\alpha\in \Theta$. 
From that description one can in particular see that each outer cone is isomorphic 
(as a cone in Euclidean space) to the face $\overline{\frak{a}^+}(\Theta) 
:= \{ H\in \overline{\frak{a}^+}\mid \beta(H)=0,\ \ \ \forall \ \beta \in _{\Bbb 
Q}\!\Delta - \Theta \}$ in the ``model chamber'' $\overline{\frak{a}^+}\subset \frak{g}$
resp.  $\overline{A^+}x_0\subset X$.

We define also the {\it outer cone}
 $CO(p)$ at $p\in \partial V(s)$ as
 the Euclidean cone in the tangent space $T_pV$ whose vertex is $0$ and whose basis is
the (closed) spherical simplex $O(p)$.

\subhead 3.3. Outer angles and structure of $V$at infinity \endsubhead

The 
 outer angles  $O(p)$ also reflect the structure at infinity of $V= \Gamma \backslash 
X$. More precisely, we shall see below how they correspond to   the spherical simplices in a 
geometric realization of the quotient modulo $\Gamma$ of the  Tits building ${\Cal T}$
 over $\Bbb Q$ associated 
to $\bold G$.

In order to describe that quotient more explicitly  we first note that the
 discrete, arithmetic group $\Gamma$ acts on $|\Cal T|$.  The resulting quotient space 
$\Gamma\backslash |\Cal T|$  
turns out to be  a finite simplicial complex. In fact, as mentioned in  Section 2,
 the conjugacy classes 
of elements of the collection $\Cal P$ of parabolic $\Bbb Q$--subgroups of $\bold G$ are in 
one-to-one correspondence with
the subsets $\Theta$ of the set $_{\Bbb Q}\Delta$
of simple $\Bbb  Q$--roots. Every conjugacy class has a standard
representative denoted by $\bold P_{\Theta}$.
 One can show that the sets  
 of double cosets $\Gamma\backslash \bold G(\Bbb Q)/
\bold P_{\Theta}(\Bbb Q)$  are
 {\it finite} for all $\Theta$ (see [B1] 15.6).
Let $\bigtriangleup \subset {\Bbb R}^q$ be the spherical $q-1$ simplex obtained by taking 
unit vectors in $\overline{\frak a^+}\subset {\frak a}\cong {\Bbb R}^q$ ($q = \Bbb Q$--rank 
of $\bold G$); note that with our notation above $\bigtriangleup\cong(\overline{A^+} 
x_0)(\infty)$. 
For  a subset $\Theta$ of  $_{\Bbb Q}\Delta$   we define the boundary simplex 
$\bigtriangleup (\Theta)$
 of $\bigtriangleup$ 
as $\bigtriangleup (\Theta) := \bigtriangleup \cap \{\alpha = 0\mid \alpha \in _{\Bbb 
Q}\!\Delta
- \Theta \}$. A vertex of $\bigtriangleup$ then corresponds to a {\it maximal}
parabolic $\Bbb Q$--subgroup while the entire simplex  
$\bigtriangleup$ corresponds to the {\it minimal} parabolic $\Bbb Q$--subgroup $\bold P := 
\bold P_{\emptyset}$ of $\bold G$.
 Let the set  $ \Gamma\backslash {\bold G}
(\Bbb  Q)/\bold P(\Bbb  Q)$ be represented by
 $\{q_1,\ldots q_m\}$ (see Proposition 2.2) and  take $m$ isometric  copies
$\bigtriangleup ^j$ of $\bigtriangleup$ 
with faces $\bigtriangleup^j(\Theta)$ corresponding to $\Theta$.
The corresponding isometries $\bigtriangleup \simeq \bigtriangleup ^j$
are denoted by $\varphi_j$.
The simplicial complex $\Gamma\backslash |\Cal T|$,
which provides a geometric realization of the 
quotient of the Tits building
$|\Cal T|$ modulo $\Gamma$,
 is constructed from
the simplices $\bigtriangleup^1, \ldots, \bigtriangleup^m$
through the following {\it incidence relations}:

Two simplices $\bigtriangleup^j$ and $\bigtriangleup^l$ are pasted
together along the faces $\bigtriangleup^j(\Theta)$ and 
$\bigtriangleup ^l (\Theta)$ by the isometry
 $\varphi_j\circ\varphi_l^{-1}\mid_{\bigtriangleup^l(\Theta)}$ if and only if 
$\Gamma q_j\bold P_{\Theta}(\Bbb Q) = \Gamma q_l \bold P_{\Theta}(\Bbb Q)$.
The simplicial complex $|\Cal T|$ and hence also $\Gamma\backslash |\Cal T|$ is connected
if the $\Bbb Q$--rank $q$ is greater than or equal to two; moreover $\dim \Gamma\backslash 
|\Cal T| = \roman{rank}_{\Bbb Q}\ \Gamma -1$.

Having constructed the finite simplical complex $\Gamma\backslash |\Cal T|$ we briefly 
explain how its simplices correspond to outer angles of the polyhedra $V(s)$ for an   
exhaustion $V = \bigcup _{s\geq 0}V(s)$ of $V$. Let ${\Cal E}$ be the set of double cosets 
 $\Gamma\backslash \bold G(\Bbb Q)/
\bold P_{\Theta}(\Bbb Q), \Theta\subset _{\Bbb Q}\Delta$. By [B1] 15.6, ${\Cal E}$ is a finite set.

By definition the elements of ${\Cal E}$  are in one to one correspondence to
the $\Gamma$--equivalence classes of parabolic ${\Bbb Q}$--subgroups of $\bold G$ or 
the simplices of $\Gamma\backslash |{\Cal T}|$ (see the above construction or [L3], Section 4). 
The elements of ${\Cal E}$  also index  the boundary strata of $V(s)$. 
In fact, if $p$ is some point in the interior of a boundary stratum of $V(s)$, then there is a 
unique set $\Theta$ of simple roots such that  the 
outer angle at  $p$ is  spanned by the gradients of the intersecting 
horospheres $ -\roman{grad}\ \tilde{h}_{i\alpha}, \alpha\in \Theta$, and hence 
$O(p)\cong \bigtriangleup(\Theta)$ (and $CO(p)\cong\overline{\frak{a}^+}(\Theta)$).

In [L3], Section 4, we showed that every point in the complement 
of the interior $V(s)^o$ of $V(s)$ in $V$ is contained in the image under the Riemannian 
exponential
map of the outer cone $CO(p)$  at some point of $\partial V(s)$. More precisely we have

\proclaim{Proposition 3.2} Let $V = \bigcup _{s\geq 0}V(s)$ be a polyhedral  
exhaustion 
as described above. Then there is 
a disjoint union
$$V -  V(s)^0 = \coprod_{p\in\partial V(s)} \roman{Exp}_p CO(p).
$$
Moreover, the points at infinity of those
unit speed geodesic rays whose initial vectors are in $O(p)$ span a (closed) spherical simplex 
$\sigma_E$, $ E\in {\Cal E}$, of the finite complex $\Gamma\backslash |{\Cal T}|$. 
\endproclaim

\subhead 3.4. The fine structure of the boundary strata\endsubhead

We will need some further details about the structure of the boundary strata of a 
polyhedral 
exhaustion.  For proofs we refer the reader to [L3], Section 4. Pick $E=E(\Theta)\in {\Cal E}$
and let $\sigma_E$ be the corresponding simplex of the finite complex $\Gamma\backslash |{\Cal T}|$. 
The pointwise isotropy group of $\sigma_E$ is isomorphic to the group of real points, say 
$P_E$,  of a conjugate under $\Gamma$ of the standard parabolic subgroup
 ${\bold P}_E$ of $\bold G$ (see [L2], Lemma 1.2).
 Recall from Section 2 that there is a Langlands decomposition 
$P_E = U_EM_EA_E$ and an associated decomposition  $X\cong U_E\times Z_E\times A_E$.
The arithmetic group $\Gamma\cap P_E$ has trivial $A_E$ factor and has a finite 
index subgroup which is the semidirect product of $\Gamma\cap U_E$
with $\Gamma\cap M_E$. 
A{ \it truncated locally symmetric space} $\Gamma\cap M_E\backslash Z_E(s)$ is a compact submanifold 
with corners in a totally geodesic subspace of $V$ isometric to $\Gamma\cap M_E\backslash Z_E$. 

\proclaim{Proposition 3.3} Let $V = \bigcup _{s\geq 0}V(s)$ be a polyhedral  
exhaustion of $V$. Then, for each $s$, the boundary $\partial V(s)$ consists of  a finite 
number of
 strata $V_E(s),\ E=E(\Theta)\in {\Cal E},$  each of which is
 a fibre bundle over a truncated
locally symmetric space with a compact nilmanifold as fibre: 
$$0\rightarrow \Gamma\cap U_E\backslash U_E\rightarrow  V_E(s)
\rightarrow \Gamma\cap M_E\backslash Z_E(s)\rightarrow 0.$$
The combinatorics of these strata is ``dual'' to that  of the  simplicial complex 
$\Gamma\backslash |{\Cal T}|$, i.e., $(n-k)$-dimensional 
 strata of $\partial V(s)$ correspond to $(k-1)$-simplices  of $\Gamma\backslash |{\Cal
T}|$.

Moreover, the outer cones $CO(p)$ at the interior points of $V_E(s)$ are all isomorphic
to some face $\overline{\frak{a}^+}(\Theta)$ (uniquely determined by the index $E=E(\Theta)$) 
of the ${\Bbb Q}$--Weyl chamber $\overline{\frak{a}^+}$.
\endproclaim

We now consider the universal covering space $X(0)\subset X$ of the polyhedron $V(0)$ of some
 exhaustion of $V$.  
The distance function $d_V$ on $V$ is defined  by
$d_V(\pi(x), \pi(y)) = \inf_{\gamma\in\Gamma} d_X(x, \gamma\cdot y)$, where $\pi : X\rightarrow 
V = \Gamma\backslash X$  denotes the canonical projection.
The next lemma will  be used in the proof of Theorem 4.1.

\proclaim{Lemma 3.4} Let
$CO(z)\subset X$ be an  outer cone of the  submanifold
 with corners $X(0)\subset X$. For $x\in CO(z)$ the geodesic segment $[\pi(x)\pi(z)]$ 
realizes the shortest distance (in $V$) from $\pi(x)$ to the compact manifold with corners
$V(0)=\Gamma\backslash X(0)$; and, in particular, $d_V(\pi(x),\pi(z))=d_X(x,z)$.
\endproclaim

\demo{Proof} Assume that there is a point $\pi(y)\in V(0)$ such that 
$$
d_V(\pi(x),\pi(y))<d_V(\pi(x),\pi(z))=d_X(x,z);
$$
where the equality comes from the fact 
(see [L3], Corollary 3.5.) that $\pi(x)$ lies on a ray in $V$ starting at $\pi(z)$.
Since $d_V(\pi(x),\pi(y))=\inf_{\gamma\in\Gamma}d_X(x,\gamma y)$ the assumption implies that 
for some $\gamma\in\Gamma$, $ d_X(x,\gamma y)< d_X(x,z)$.  Since $X(0)$ is 
$\Gamma$-invariant and $y\in X(0)$ we have  $\gamma y\in X(0)$ and hence $d_X(x,X(0))<d_X(x,z)$.
This is a contradiction since the segment $[xz]$ clearly realizes the shortest distance in $X$ from
$x$ to $X(0)$.
\qed\enddemo

We conclude this section with an interesting  metric property of 
$\Gamma\backslash |\Cal T|$ which will be used in Section 5. The definition of a ``CAT($1$) 
space'' can be found for example in [BrHa] or [Ba].

\proclaim{Lemma 3.5} The finite spherical simplical complex
 $(\Gamma\backslash |\Cal T|, d_{\Cal T})$
is a $\roman{CAT}(1)$ space. 
\endproclaim

\demo{Proof} The model simplex $\bigtriangleup$  and all its faces are  convex subsets of 
the unit sphere
in ${\Bbb R}^q$
and therefore  CAT(1) spaces. Pasting CAT(1) spaces along complete  $(\pi)$--convex 
subspaces yields again CAT(1) spaces (see e.g. [BrHa], Ch.II. 4.4).
\qed\enddemo

\head 4. Metric properties of  arithmetic quotients of symmetric spaces \endhead

\subhead 4.1. Isometric images of Weyl chambers\endsubhead

In this section we derive estimates which are crucial for the rest of the paper. 
We will prove that the canonical projection map $\pi:X\longrightarrow V$ is an isometry 
(with repect to the distance function $d_X$ and $d_V$) when restricted to 
certain  closed ${\Bbb Q}$--Weyl chambers.

\proclaim{Theorem 4.1} Let $V = \bigcup _{s\geq 0}V(s)$ be an exhaustion of $V$ by manifolds 
with 
corners (as in Section 3)
and let $CO(p), p\in V(0)$, be any outer cone of $V(0)$.
Then  the restriction of the canonical projection $\pi:X\longrightarrow V$  
to any lifted outer cone $CO(z)\subset X$, with $\pi(z)=p$, is an isometry with respect
 to the distance functions of $X$ resp. $V$.

 In particular, there is a $\Bbb Q$--Weyl chamber ${\Cal C}_*\subset {\Cal C}$ 
such that for every $i\in\{1,\ldots,m\}$ the chamber $q_i{\Cal C}_*\subset \Omega\subset X$ is 
isometrically mapped to $V$ under  $\pi$.
\endproclaim

\demo{Proof} We consider an arbitrary exhaustion $V = \bigcup _{s\geq 0}V(s)$.
For a point $\pi(z)\in V(0)$ let $CO(\pi(z))$ be the outer cone. In the lift 
of $CO(z)$ of $CO(\pi((z))$ to $X$  we pick two arbitrary points $x,y$.
 We wish to show that $d_V(\pi(x),\pi(y))\geq d_X(x,y)$ (the opposite inequality holds
by definition). Let $c$ be a  curve of minimal length in $V$ between $\pi(x) $ and $\pi(y)$, i.e.,
 $L(c)=d_V(\pi(x),\pi(y))$.
We distinguish two cases: (1) $c$ intersects the compact polyhedron 
$V(0)$ and (2) $c$ does not intersect $V(0)$.

Case (1): Decompose $c$ into three segments $c=c_1\cup c_2\cup c_3$, such that $c_1$ connects 
$\pi(x)$ to the first
intersection point of $c$ with $V(0)$ and  $c_3$ connects the last intersection point 
of $c$ with $V(0)$ to $\pi(y)$. By Lemma 3.4. we have 
$$
d_V(\pi(x), V(0))=d_V(\pi(x),\pi(z))=d_X(x,z)\  \  \roman{and}$$
 $$ d_V(\pi(y), V(0))=d_V(\pi(y),\pi(z))=d_X(y,z).
$$
 Hence
$$
d_V(\pi(x),\pi(y))=L(c)=L(c_1)+L(c_2)+L(c_3)\geq d_X(x,z)+d_X(y,z)\geq d_X(x,y).
$$

Case (2): According to Propositions 3.2 and 3.3  we can decompose the curve 
$$
c: [0,d_V(\pi(x),\pi(y))]\longrightarrow V\setminus V(0)
$$
 into a finite
number of segments $c_i, 1\leq i\leq l$, say, each of which is contained in a subset
of $V$ diffeomorphic to  $W_E(0)\times \overline{\frak{a}^+}(\Theta)$, where $W_E(0)$ is a nilmanifold
bundle over truncated locally symmetric and compact and the ${\Bbb Q}$-Weyl chamber face 
$\overline{\frak{a}^+}(\Theta)\cong \overline{A^+}(\Theta)x_0$ 
corresponds to the (congruent) outer cones with $E=E(\Theta)$. The Riemannian metric  
 on $X$ with respect to 
horocyclic coordinates corresponding to  $A_E$ supporting $\overline{A^+}(\Theta)$ is given by  
a  formula of Borel (see [B2], Proposition 1.6). From the $\Gamma$-invariance of the exhaustion in $X$ 
(compare the proof of Theorem 4.2 in [L2])
we deduce that the Riemannian metric on $V$ at the point 
$(w,a)\in W_E(0)\times \overline{A^+}(\Theta)$  can be written in the form
$$
ds_{(w,a)}^2=dw_{(w,a)}^2+da_a^2.
$$
For the length of the segment $c_i(s)=(w_i(s),a_i(s)), s\in [s_i,s_{i+1}]$, we thus have
the estimate $L(c_i)\geq L(a_i), 1\leq i\leq l$. Hence
$$
d_V(\pi(x),\pi(y))=L(c)=\sum_{i=1}^l L(c_i)\geq  \sum_{i=1}^l L(a_i).
$$
Now consider the continuous curve $a(s):=\bigcup_{i=1}^l a_i(s)$ in the chamber 
$\overline{A^+}\subset A$. If we write $x=ax_0$ and $y=bx_0$, the curve $a(s)x_0$ connects
$x$ and $y$ in $\overline{A^+}x_0$ and satisfies
$$
\sum_{i=1}^l L(a_i)\geq d_A(a,b)=d_X(x,y).
$$
Together with the estimate above this proves the claim also in case (2).

 Finally, the chamber  ${\Cal C}_*$ can be chosen as follows. Consider the preimages
 in the fundamental set $\Omega$ 
 of the minimal boundary strata of $\partial 
V(0)$. They  intersect the $\Bbb Q$--Weyl chambers $q_i {\Cal C}$, $i=1,\ldots,m$, in well defined  
points $x_i$ (see [L2] 
Lemma 3.5). The lifted  (maximal) outer cones  based at the points $\pi(x_i)$ of $\partial 
V(0)$ are again $\Bbb Q$--Weyl chambers $q_i{\Cal C}_*\subseteq q_i{\Cal C}$, $i=1,\ldots, 
m$, 
and  in finite Hausdorff-distance from $q_i{\Cal C}$. 
\qed \enddemo

\subhead 4.2.  A net in $V$\endsubhead

Theorem 4.1 yields a metric space which is quasi-isometric to $V$. Recall that a
 subset $N$ of a metric 
space $(S,d)$ is called a ($\varepsilon$--){\it net} if there is some positive  constant 
$\varepsilon$ such that  $d(s,N)\leq \varepsilon$ for all $s\in S$; in particular the 
Hausdorff-distance
between $N$ and $S$ is at most $\varepsilon$.
 
\proclaim{Corollary 4.2} There is a net $N$  in $V$ consisting of finitely  many
 isometrically embedded  $\Bbb Q$--Weyl chambers.
\endproclaim

\demo{Proof} Let ${\Cal C}_*$ be as in Theorem 4.1. The definition of a fundamental set as a 
finite  union of Siegel sets
$\Omega = \bigcup_{i = 1}^m q_i\cdot{\Cal S}$ (see Proposition 2.2) implies that $\bigcup_{i 
= 1}^m q_i\cdot{\Cal C}$
and hence also $\bigcup_{i = 1}^m q_i\cdot{\Cal C}_*$ is a net in $\Omega\subset X$. Since 
 $V = \pi(\Omega)$ and $d_V(\pi(x),\pi(y))\leq d_X(x,y)$ for all $x,y\in X$ the claim 
follows.
\qed\enddemo

\remark{\bf Remark 4.3} From Proposition 3.1 (see also [L2], Lemma 2.4, 2.5)
it follows that the $q_i{\Cal C}_*$, 
$1\leq i\leq m$, can be chosen to be  disjoint in $X$ and such that $\pi: X\rightarrow V$ 
restricted to $\bigcup_{i=1}^m q_i{\Cal C}_*$ is bijective.
\endremark

\head 5. The asymptotic cone  and  the proof of a conjecture of 
Siegel \endhead

\subhead 5.1. A quasi-isometric model for a locally symmetric space\endsubhead

The goal of  this section is to show that  the asymptotic cone of a non-compact locally 
symmetric
space $V = \Gamma\backslash X$ of finite volume is isometric to the Euclidean cone over a 
finite simplicial complex, the quotient of the Tits bulding $\Cal T$ of ${\bold G}$ modulo 
$\Gamma$.

If $\text{rank}_{ \Bbb Q }\Gamma \geq 2$ the geometric realization  $|\Cal T|$  of $\Cal T$ 
as a spherical complex carries a natural distance function 
which we normalize such that the diameter of $|\Cal T|$ is $\pi$.  This ``spherical'' metric 
 induces a distance function    $d_{\Cal T}$ on the 
finite simplical complex $\Gamma\backslash |\Cal T|$.  
We define the {\it Euclidean cone}  $C(\Gamma\backslash |\Cal T|)$  {\it over }
 $\Gamma\backslash |\Cal T|$ to be the product $[0,\infty)\times  \Gamma\backslash |\Cal T|$
with $\{0\}\times \Gamma\backslash |\Cal T|$ collapsed to a point $\Cal O$  and  endowed 
with the
 cone metric
$$d_C^2((a,x),(b,y)) := a^2 + b^2 - 2ab\cos d_{\Cal T}(x,y).$$
If  $\text{rank}_{ \Bbb Q }\Gamma =1$ we set $d_C((a,x),(b,y)) := a+b$ if $x\neq y$ and 
$d_C((a,x),(b,y)):=|a-b|$ if $x=y$.

In the same way one defines the Euclidean cone over an arbitrary metric space
 of diameter $\leq \pi$ (see e.g. [BrHa]).

There is an alternative way to build up $(C(\Gamma\backslash |\Cal T|), d_C)$ which
parallels the construction of $\Gamma\backslash |\Cal T|$ in Section 3:
We there realized the spherical simplices 
 $\bigtriangleup^j$ for $j= 1,\ldots,m$  on the unit sphere in ${\Bbb R}^q$; thus we can 
take the cones $C(\bigtriangleup^j)$ in the Euclidean space ${\Bbb R}^q$ with vertex $0$. We 
endow these simplicial cones with the induced Euclidean metric and glue them together using 
the same
incidence relations as for $\Gamma\backslash |\Cal T|$.

We recall the notion of Hausdorff-convergence of (unbounded) pointed metric spaces
 (see [GLP] or [Bu]).
The {\it distortion} of a map $f: A\rightarrow B$ of metric spaces $A$ and $B$  is defined 
as
$$\roman{dis}(f) := \sup_{a,b}|d_A(a,b) - d_B(f(a), f(b))|.
$$
 The {\it uniform distance}
between metric spaces $A$ and $B$ is defined as $|A,B|_u = \inf_f\roman{dis}(f)$ where the 
infimum is taken over all bijections $f:A\rightarrow B$.
A sequence of metric spaces $M_n$ {\it Hausdorff-converges} to a metric space $M$ iff for 
every $\varepsilon>0$ there is an $\varepsilon$--net $M_{\varepsilon}$ in $M$ which is the 
uniform limit of $\varepsilon$--nets $(M_n)_{\varepsilon}$ in $M_n$.
We say that a sequence $(M_n,p_n)$ of unbounded, pointed metric spaces Hausdorff-converges 
to a pointed 
metric space $(M,p)$ if for every $r>0 $ the balls $B_r(p_n)$ in $M_n$  Hausdorff-converge 
to the ball $B_r(p)$ in $M$.

Let $v_0$ be an (arbitrary)  point of the locally symmetric space $V=\Gamma\backslash X$.
We define the {\it asymptotic cone} of $V$ as the Hausdorff-limit of pointed metric spaces:
$$ \text{Cone}(V) := {\Cal H}-{\lim}_{n\rightarrow\infty} (V, v_0, \frac{1}{n} d_V).$$

By Corollary 4.2 there is a net $N$ in $V$ consisting of disjoint closed Euclidean cones
$\pi(q_i{\Cal C}_*))$ in $V$. We identify
the abelian Lie algebra ${\frak a}$ with ${\Bbb R}^q$. For all $n\in \Bbb N$ we define a map 
$$f_n: N\subset (V, \frac{1}{n}d_V)\rightarrow 
C(\Gamma\backslash |\Cal T|);\ \text{by}$$
$$ f_n(\pi(q_j(\exp H)\cdot x_0)) := \frac{1}{n}H \in 
 \log {\Cal C}_* \cong C(\bigtriangleup^j)\subset {\Bbb R}^q.
$$
 By Remark 4.3 and the second construction of $C(\Gamma\backslash |\Cal T|)$  $f_n$ is a 
bijection from the interior of $N$ onto its image  in $C(\Gamma\backslash |\Cal T|)$;
note that this image  is open and dense in $C(\Gamma\backslash |\Cal T|)$ for all $n$.

\proclaim{Lemma 5.1} There is a constant $D\geq 0$ such that 
$\roman{dis}(f_n)\leq \frac{1}{n}D$ for all $n\in \Bbb N$, i.e., $f_n$ is a 
$(1,\frac{1}{n}D)$--quasi-isometry. In particular,  $V$ is quasi-isometric to   
$C(\Gamma\backslash |\Cal T|)$.
\endproclaim

\demo{Proof} We consider some polyhedral exhaustion $V = \bigcup_{s\geq 0}V(s)$. The 
intersection of the net $N$ with  $V\setminus V(0)$ is still a net in $V$ since $V(0)$ is 
compact. Let $u,v$ be two points  in the interior of $N\cap (V\setminus V(0))$.
We take a path $c([0,L])$ in $V$ between $u$ and $v$ of minimal length and parametrized by 
arc-length. By Proposition 2.2  we can write  $V = \pi(\Omega) = \bigcup_{j = 1}^m q_j{\Cal 
S}$ for a fundamental set $\Omega$. We can thus find $t_i\in[0,L]$, $ 0\leq i\leq n+1$ such 
that
$t_0 = 0, t_{n+1}= L$ and   $[0,L] = \bigcup_{i=0}^n[t_i,t_{i+1}]$ with $c([t_i, 
t_{i+1}])\subset \pi(q_{j_i}{\Cal S})$ for $j_i\in\{1,\ldots,m\}$. Associated to the path 
$c$ we thus get a well-defined string
$(j_0,j_1,\ldots,j_n)$. 
We next replace $c$ by a path $\overline{c}$ whose associated string contains an element 
$j_k\in  \{1,\ldots,m\}$ at most once: Start with $j_0$ and assume that $0\leq l\leq n$ is 
the greatest index
such that $j_l = j_0$. Geometrically this means that the path $c$ returns to 
$\pi(q_{j_0}{\Cal S})$ at $c(t_l)$. Using   Theorem 4.1 and the 
definition of a Siegel set $\Cal 
S$ we can replace $c([t_0,t_{l+1}])$ by a geodesic segment, say 
$\overline{c}([s_{j_0},s_{j_0+1}])$, in $\pi(q_{j_0}{\Cal C})$ whose length $L$  satisfies
$$L(\overline{c}([s_{j_0},s_{j_0+1}]))\leq L(c[t_0,t_{l+1}]) + 2D_1,$$
for  $D_1:=\text{diam}(\omega)$. Repeating this procedure with $j_{l+2}$ etcetera we eventually 
get 
a sequence of at most $m$ segments 
$\overline{c}([s_k,s_{k+1}])$ in $N$ of total length $\leq L(c) + 2mD_1.$
We next show that these segments can be chosen in such a way that they are mapped by $f_1$ to a 
continous path 
$\tilde{c}$ in $C(\Gamma\backslash |\Cal T|)$ from $f_1(u)$ to  $f_1(v)$.
In fact,  each $\pi(q_j{\Cal C}_*)\subset \pi(q_j{\Cal S})\subset V$ is isometric to a
simplicial cone in $C(\Gamma\backslash |\Cal T|)$. The map $f_1$ thus yields  well defined 
images of (the interior of)  the segments
$\overline{c}([s_k,s_{k+1}])\in \pi(q_k{\Cal C}_*)$. By construction the endpoint 
$\overline{c}(s_{k-1})\in \pi(q_{k-1}{\Cal C}_*)$ and the initial point 
$\overline{c}(s_{k})\in \pi(q_{k}{\Cal 
C}_*)$ are on the same levelset, say $\partial V(s)$, of the exhaustion
 of $V$ we have chosen above and by construction are at most the  distance $2D_1$ apart.
By the (second) construction of $C(\Gamma\backslash |\Cal T|)$ we can join their  images
 by segments in $C(\Gamma\backslash |\Cal T|)$ of uniformly bounded length 
 $2D_1$ to obtain the path $\tilde{c}$. This argument yields the 
 estimate
$$
d_C(f_1(u),f_1(v)) \leq L_C(\tilde{c})\leq 4mD_1+ L(c)= 4mD_1 + d_V(u,v).
$$

On the other hand, given a geodesic path  $\tilde{c}$ in $C(\Gamma\backslash |\Cal T|)$ 
joining
two  points $x$ and $y$ in the image of the interior of  $N\cap (V\setminus V(0))$ under  
$f_1$, we can lift the segments contained in the simplicial cones 
$C(\bigtriangleup^j)$ to the corresponding chambers $\pi(q_j{\Cal C}_*)$ via $f_1^{-1}$. By 
the same 
arguments as before the endpoints of the lifted segments can be joined in $V$ to form a 
continuous path $c$ between $f_1^{-1}(x)$ and $f_1^{-1}(y)$ of length
$$
d_V(f_1^{-1}(x),f_1^{-1}(y)) \leq L(c) \leq  2mD_1 + L(\tilde{c}) \leq 2mD_1 + d_C(x,y).
\tag 2$$
Combining $(1)$ and $(2)$ and setting $D:= 4mD_1$ we get for all $u,v$ in (the interior of)
 $N\cap(V\setminus V(0))$
that
$$|d_V(u,v) - d_C(f_1(u), f_1(v))|\leq D.$$
Finally since for any $n\in\Bbb N$ 
$$\split |\frac{1}{n}d_{V}(u,v) - d_C(f_n(u), f_n(v))| &= |\frac{1}{n}d_{V}(u,v) - 
\frac{1}{n}d_C(f_1(u), f_1(v))| = \\
 &= \frac{1}{n} |d_V(u,v) - d_C(f_1(u), f_1(v))| \leq \frac{1}{n}D,
\endsplit
$$
  the Lemma follows.
\qed\enddemo

We are now prepared to state a key result of this paper. The following theorem was also 
proved
independently by T. Hattori using  equivariant Hausdorff limits with respect to suitable 
finite index subgroups of $\Gamma$ and an embedding into a space of positive definite matrices
(see [H3]).
 
\proclaim{Theorem 5.2}  Let $X$ be a Riemannian symmetric space of rank $\geq 2$, let 
$\Gamma$ be an irreducible, non-uniform  lattice in the group of isometries of $X$  and let 
$V$ be the locally symmetric space  $\Gamma\backslash X$. 
Then the asymptotic cone $\roman{Cone}(V)$ is  isometric  to 
the Euclidean cone $C(\Gamma\backslash |\Cal T|)$ over the finite simplicial complex 
$\Gamma\backslash |\Cal T|$. 
\endproclaim

\demo{Proof} Given $\varepsilon >0$ there is $n_{\varepsilon}\in\Bbb N$, such that 
$N' := N\cap(V\setminus V(0))\subset V$ is an $\varepsilon$--net in $(V, \frac{1}{n}d_V)$ 
for all $n\geq n_{\varepsilon}$.
Since $f_1(N')$ is dense in $C(\Gamma\backslash |\Cal T|)\setminus f_1(N\cap V(0))$ it is an 
$\varepsilon$--net in the latter for all
$\varepsilon>0$. Then  there is $r_0$ such that for  $r\geq r_0$  the same assertions are 
true 
for the subsets of balls with radius $r$:
$$
N'\cap B_r(v_0)\subset (V, \frac{1}{n}d_V)\ \roman{and}\  f_n(N')\cap B_r(\Cal O)
\subset C(\Gamma\backslash |\Cal T|).
$$
From Lemma 5.1 and its proof we obtain the following uniform estimate:
$$
|N'\cap B_r(v_0),f_n(N')\cap B_r(\Cal O)|_u\leq \frac{1}{n}D,
$$
which implies the theorem.
\qed\enddemo

The definition of the asymptotic cone of a metric space $S$ implies that  bi-lipschitz 
invariants of  $\roman{ Cone}(S)$ are quasi-isometry invariants of $S$
(see [G] for a discussion). 
It is well-known that the $\Bbb R$--rank of $\bold G$ equals the rank of the symmetric space 
$X$,
i.e., the dimension of a maximal totally geodesic submanifold of $X$ with sectional 
curvature zero. The dimension of a simplicial complex is equal to the dimension of a maximal
simplex. We thus define 
 $$\dim \roman{Cone} (\Gamma\backslash X) := \dim (\Gamma\backslash |\Cal T|) +1.$$
 Theorem 5.2 yields the following geometric interpretation of the $\Bbb Q$--rank and the 
cohomological dimension, $\text{cd}\, \Gamma$, of $\Gamma$:

\proclaim{Corollary 5.3} Let $X$ be  a symmetric space of noncompact type and  rank $\geq 2$ 
and let $\Gamma$
be an irreducible (arithmetic) lattice in the isometry group of $X$. 
 Then 
$$\align
\roman{rank}_{\Bbb Q}\Gamma &= \dim \roman{Cone} (\Gamma\backslash X)\tag i\\
\roman{cd}\,\Gamma &= \dim X- \dim \roman{Cone} (\Gamma\backslash X).\tag ii
\endalign
$$
\endproclaim

\demo{Proof}  By a well-known result of Borel and Harish-Chandra $V$ is compact if and only 
if $\roman{rank}_{\Bbb Q}\Gamma = 0$ (see [BH]). If $V = \Gamma\backslash X$ is compact then 
$\roman{Cone} (V)$ is a point, i.e., 
$\dim \roman{Cone} (V) = 0$. On the other hand if $V$ is not compact there is a geodesic ray
in $V$ and thus  $\dim \roman{Cone} (V) \geq 1$. This proves the formula in the compact 
case. If $V$ is not compact, i.e., if  
 $\roman{rank}_{\Bbb Q}\Gamma \geq 1$, we have by
 Theorem 5.2 
$$\dim \roman{Cone} (\Gamma\backslash X) = \dim (\Gamma\backslash |\Cal T|) +1
= \dim |\Cal T| + 1 = \roman{rank}_{\Bbb Q}\bold G = \roman{rank}_{\Bbb Q}\Gamma
$$
and hence (i). Equality (ii) follows from (i) and [BS] 11.4.3.
\qed\enddemo

\remark{\bf Remark 5.4} In [KL]  B. Kleiner and B. Leeb showed that the asymptotic cone
(with respect to any fixed ultrafilter) of a globally symmetric space $X$  is a Euclidean 
building. The 
rank of this building is equal to the rank of $X$, i.e., to the $\Bbb R$--rank of $\bold G$. 
The latter is therefore a quasi-isometry invariant of $X$, whereas the $\Bbb Q$--rank is a 
quasi-isometry invariant of $\Gamma\backslash X$.
\endremark

\vskip 1pt

\proclaim{\bf Corollary 5.5} The asymptotic cone of a locally symmetric space
$V=\Gamma\backslash X$ is a $\roman{CAT}(0)$ space and thus in particular contractible.
\endproclaim

\demo{Proof} By Theorem 5.2 $\roman{Cone}(V)$ is isometric to the Euclidean cone 
$C(\Gamma\backslash |\Cal T|)$. A theorem of V.N. Berestovski says that for a geodesic 
metric 
space $Y$ the Euclidean cone $C(Y)$ is a CAT(0) space  iff $Y$ is a CAT(1) space (see 
[BrHa], Ch. II.4). The claim then follows from Lemma 3.4.
\qed\enddemo

\subhead 5.2. A conjecture of Siegel and its proof\endsubhead

\remark{\bf Remark 5.6}  Theorem 4.1 could be deduced from an estimate about 
the distance function  $d_X$ of $X$ when restricted to a fundamental set $\Omega$, which 
is given in a paper by A. Borel (see [B2], Theorem 2.3). But unfortunately
the second part of the proof given there only works in the $\Bbb Q$--rank 
$1$ case and is not  complete for $\Bbb Q$--rank $\geq 2$ (this was also 
observed in [JM]). More precisely,
in [B2], pp. 550-552, inequality $(12)$ does not imply inequality $(14)$ and this in
turn does not imply inequality $(5)$ as is claimed there. Nevertheless the assertion of the 
theorem, which in particular answers a question of C.L. Siegel (see [Si], Section 10), is 
correct and 
we next  present a complete proof based on Theorem 4.1 and Lemma 5.1. 
\endremark

\topinsert
$$
\psset{unit=0.7cm}
\pspicture(-7,-6.6)(3.8,6.6)
\pspolygon[linewidth=1.2pt](0,0)(3.8,5.8)(-3.3,6.6)(0,0)
\psline[linewidth=1.2pt](0,0)(-7,0)
\psline[linewidth=1.2pt](-7,0)(-3.3,6.6)
\psline[linewidth=1.2pt](-7,0)(-3.3,-6.6)
\psline[linewidth=1.2pt](0,0)(-3.3,-6.6)
\psline(-1.5,-1)(-2,0)(-1.75,3.5)(1.9,3.7)
\psline[linestyle=dotted,dotsep=2pt](-2,0)(-1.8,-3.6)(-4.3,-0.5)
\psdots(-1.5,-1)(-2,0)(-1.75,3.5)(1.9,3.7)(-1.8,-3.6)(-4.3,-0.5)
\rput(1.7,4){$q_jx$}
\rput(0,4){$3$}
\rput(-2.3,2){$2$}
\rput(-1,-.5){$1=1'$}
\rput(-1.2,-1.5){$q_iy$}
\rput(-2.3,-1.7){$2'$}
\rput(-3.6,-2){$3'$}
\rput(-4.9,-0.8){$q_ix$}
\endpspicture
$$
\botcaption{Figure  5.1}\endcaption
\endinsert

\proclaim{Theorem 5.7}  For  $q_k\in {\bold G}(\Bbb Q)$ ($k=1,\ldots,m$) and a Siegel set 
${\Cal S}$  as in Proposition 2.2 there is a
constant $D'$ such that
$$
d_X(\gamma q_ix, q_jy)\geq d_X(x,y) + D',
$$ 
for all $i,j\in\{1,\ldots,m\}$, all $\gamma\in \Gamma$ and all $x,y \in {\Cal S}$.
In particular the canonical projection $\pi: X\longrightarrow \Gamma\backslash X$ 
restricted to each Siegel set $q_i{\Cal S}$ is a $(1,D´)$--quasiisometry.
\endproclaim

\demo{Proof} The definition of the Siegel set ${\Cal S}$ implies that 
 the $\Bbb Q$--Weyl chamber ${\Cal C}_*$ is a net in ${\Cal S}$. Hence it suffices to prove 
the estimate for $x,y\in {\Cal C}_*$. By definition of $d_V$ and Lemma 5.1 we have
$$
\split 
d_X(\gamma q_ix, q_jy) &\geq d_V(\pi(q_ix), \pi(q_jy))\geq \\
                     &\geq d_C(f_1(\pi(q_ix)), f_1(\pi(q_jy))) - D.
\endsplit
$$
We can further decompose a shortest path $c$ in  $C(\Gamma\backslash |\Cal T|)$ into 
segments contained 
in the simplicial cones $C(\bigtriangleup^l)\subset  C(\Gamma\backslash |\Cal T|)$. Since 
all these Euclidean cones are isometric (for $l=1,\ldots,m$) we can take copies of 
successive  segments 
to build a continuous path $c'$ of the same length as $c$ in a single cone, say 
$C(\bigtriangleup^i)$. In  Fig. 5.1 is indicated the ``shadow'' of $c$ and $c'$ in the basis 
$\Gamma\backslash |\Cal T|$ of the cone $C(\Gamma\backslash |\Cal T|)$.
It then follows that
$$
d_C(f_1(\pi(q_jx)), f_1(\pi(q_iy)))\geq d_C(f_1(\pi(q_ix)), f_1(\pi(q_iy)))\geq 
d_V(\pi(q_ix),\pi(q_iy)) - D,
$$
where for the last inequality we used Lemma 5.1 again. Finally, by Theorem 4.1, we have that
$d_V(\pi(q_ix),\pi(q_iy)) = d_X(q_ix,q_iy) = d_X(x,y)$; this proves the theorem.
\qed\enddemo

Theorem 5.7 has independently  been proved by Ji and in a special case by Ding (see Section 1).

\head 6.  On the Tits geometry  of  locally symmetric spaces  \endhead

We consider a locally symmetric space $V = \Gamma\backslash X$ as in the previous 
section 5, i.e., $X$ is a symmetric space of  non-compact type and rank $\geq 2$ and 
$\Gamma$ is an irreducible, non-uniform (arithmetic) lattice in the isometry group of $X$. 
We remark in passing that the boundary at infinity  $\partial _{\infty}X$ of the {\it 
globally} symmetric space $X$  can be identified with
a geometric realization of 
the spherical Tits building associated to ${\bold G}(\Bbb R)$; its dimension is equal
to $\roman{rank}\, X - 1 = \roman{rank}_{\Bbb R} {\bold G} - 1$.
The following combinatorial description of the boundary at infinity of the {\it locally} 
symmetric space $V$ in terms of the Tits building 
$\Cal T$ associated to  ${\bold G}(\Bbb Q)$ is due to Ji and MacPherson (see [JM]).
As an application of Theorem 5.2  we obtain  a new proof of this result.

\proclaim{Theorem 6.1} 
The equivalence classes of asymptotic  geodesic rays in the locally symmetric space $V$ 
correspond bijectively to the points  of $\Gamma\backslash |\Cal T|$, i.e., 
$$\partial_{\infty}V\cong\Gamma\backslash |\Cal T|.$$
\endproclaim

Before proving this theorem we describe the geodesic rays in the Euclidean cone 
$C(\Gamma\backslash |\Cal T|)$.

\proclaim{Lemma 6.2} Any geodesic ray $c$ in 
$C(\Gamma\backslash |\Cal T|)$ is of the form $c(t) = (t_0+t,z)$
for some $t_0\geq 0$ and $z\in \Gamma\backslash |\Cal T|$. Moreover two rays $c_1$ and $c_2$ 
are asymptotic iff there are $t_1,t_2\geq 0$ with $c_1(t)=(t_1+t,z)$ and $c_2(t)=(t_2+t,z)$, 
i.e., $\partial_{\infty} C(\Gamma\backslash |\Cal T|)\cong \Gamma\backslash |\Cal T|$. 
\endproclaim

\demo{Proof} Let  $c: [0,\infty)\rightarrow C(\Gamma\backslash |\Cal T|)$ be a geodesic ray, 
i.e., a curve which realizes the distance between any two of its points. If $c$ eventually 
lies in a single chamber
$C(\bigtriangleup^k)$ considered as a subset of $C(\Gamma\backslash |\Cal T|)$, then it is 
clearly of the claimed form. We may  thus assume, that the ray $c$  does not
eventually stay in a  single chamber. It then has to return to a fixed chamber, say 
 $C(\bigtriangleup^j)$, i.e., there are $0\leq t'<t''$  and $1\leq j\leq m$ such that
$c(t'), c(t'')\in  C(\bigtriangleup^j)$.
Since the simplicial cones $C(\bigtriangleup^k)\subset  C(\Gamma\backslash |\Cal T|)$, 
$k=1,\ldots,m$, are all isometric  we can take  copies of the segments of $c([t',t''])$
in different chambers to build a continuous path of the same length $d_C(c(t'),c(t'')) = 
t''-t'$ in  $C(\bigtriangleup^j)$ (compare the proof of Theorem 5.7). This is a 
contradiction, because by construction (and the assumption about $c$) this path is strictly 
longer than the straight line in $C(\bigtriangleup^j)$ between the points $c(t')$ and 
$c(t'')$ on the geodesic ray $c$.
\qed\enddemo

\demo{Proof of Theorem 6.1} Let $z\in \partial_{\infty}V$ be represented by a unit-speed 
geodesic ray
 $c:[0,\infty)\rightarrow V$, i.e., $z = c(\infty)$. 
Then $c$ is also a
geodesic ray in the rescaled space $(V,\frac{1}{n}d_V)$ and hence  $c$ converges to a ray 
$\hat{c}$ in 
$\roman{Cone}(V) = \lim_{n\rightarrow\infty} (V, v_0, \frac{1}{n} d_V).$
If one also has $z=\tilde{c}(\infty)$ then the rays $\tilde{c}$  and $c$ are in finite 
Hausdorff
distance in $V$ and thus converge to the same ray $\hat{c}$ in $\roman{Cone}(V)$.
By Theorem 5.2 $\roman{Cone}(V)$ is the Euclidean cone $C(\Gamma\backslash |\Cal T|)$. 
Using Lemma 6.2 we thus get a well defined map
$$R: \partial_{\infty}V\longrightarrow  \Gamma\backslash |\Cal T|;
\ c(\infty)\longmapsto \hat{c}(\infty),
$$
which we claim to be bijective.

We first show that $R$ is {\it onto}: 
Let $z\in \Gamma\backslash |\Cal T|$. The simplicial complex
$\Gamma\backslash |\Cal T|$ was constructed by pasting together finitely many simplices 
$\bigtriangleup^k$. Thus there is at least one simplex, say $\bigtriangleup^j$, such that
$z\in  \bigtriangleup^j\subset \Gamma\backslash |\Cal T|$. We have a natural identification 
$q_j{\Cal C}_*(\infty)
\cong   \bigtriangleup^j$ (see Section 3) and by Theorem 4.1 there is a ray of the form 
$\pi(q_jc(t))\subset \pi(q_j{\Cal C}_*)\subset V$ with $R(\pi(q_jc(\infty))) = 
\hat{c}(\infty) = z$.

Next we show that $R$ is {\it one-to-one}:
Assume that  for two rays $c_1$ and $c_2$ in the net $N\subset V$ one has
$$
\hat{c}_1(\infty) = R(c_1(\infty)) = z =  R(c_2(\infty))  = \hat{c}_2(\infty).
$$
By Theorem 5.2 and Lemma 6.2  the rays $c_1$ (resp. $c_2$) Hausdorff-converge to 
$\hat{c}_1(t) = (t_1+t,z)$ (resp. $\hat{c}_2(t) = (t_2+t,z)$). This implies
that there is $n_0\in\Bbb N$ and some constant $\kappa$, such that
for all $n\geq n_0$ and all $t\in [0,1]$ one has (see proof of Lemma 5.1)
$$
\frac{1}{n}d_{V}(c_1(tn), c_2(tn))\leq  \frac{\kappa}{n}.
$$
Consequently the rays $c_1$ and $c_2$ are  asymptotic in $V$.
\qed\enddemo

Let  $X$ be a simply connected, complete Riemannian manifold of nonpositive curvature and 
let 
$z_1,z_2\in \partial_{\infty}X$ represented by geodesic rays $c_1$ and $c_2$, respectively. 
The {\it Tits metric} on $\partial_{\infty}X$ is  defined  by the formula
$$ 2\sin\frac{1}{2}d_{\Cal T}(z_1,z_2) := 
\lim_{t\rightarrow\infty}\frac{1}{t}d_X(c_1(t),c_2(t)).$$ 
This generalizes the corresponding formula for  a globally symmetric space of non-compact 
type
where $d_{\Cal T}$ {\it is} the Tits distance of the spherical Tits building of ${\bold 
G}(\Bbb R)$.
See the book [E]  for  details.
The next theorem asserts that exactly the same intrinsic formula also holds in the  case of 
non simply connected, arithmetic quotients. It confirms a conjecture of Ji and MacPherson
 (see [JM], 6.4).

\proclaim{Theorem 6.3} For two points $z_1,z_2$ in $\partial_{\infty}V\cong \Gamma\backslash 
|\Cal T|$ represented by unit 
speed geodesic rays $c_1(t), c_2(t)$ in $V$, respectively, one has the formula
$$ 2\sin\frac{1}{2}d_{\Cal T}(z_1,z_2) = 
\lim_{t\rightarrow\infty}\frac{1}{t}d_V(c_1(t),c_2(t)).$$
\endproclaim

\demo{Proof} By Theorem 5.2 and the definitions the balls of radius $1$ in 
$(V,\frac{1}{t}d_V)$
centered at $v_0$ Hausdorff-converge to the ball of radius $1$ centered at $\Cal O$ in 
$C(\Gamma\backslash |\Cal T|)$. Hence using Theorem 6.1 we have
$$
\split
\lim_{t\rightarrow \infty}\frac{1}{t}d_V(c_1(t),c_2(t)) & =  d_C((1,z_1),(1,z_2)) =\\
& =[2 - 2\cos d_{\Cal T}(z_1,z_2)]^{1/2} = 2\sin\frac{1}{2}d_{\Cal T}(z_1,z_2).
\endsplit
$$
Alternatively, one can identify the rays $c_1$, $c_2$ with rays in the cone 
$C(\Gamma\backslash |\Cal T|)$ and then use Corollary 5.5 and the fact the formula holds in 
any CAT(0) space
(see [BrHa], Ch.III 3.6).
\qed\enddemo

\enddocument